\documentclass[11pt]{article}
\usepackage[T1]{fontenc}
\usepackage{lmodern}
\usepackage[a4paper,margin=27mm]{geometry}
\usepackage{amsmath,amssymb,amsthm,mathtools}
\usepackage{microtype}
\usepackage{authblk}
\usepackage[hidelinks]{hyperref}
\usepackage{bookmark}

\newtheorem{theorem}{Theorem}[section]
\newtheorem{lemma}[theorem]{Lemma}

\newtheorem{conjecture}[theorem]{Conjecture}
\theoremstyle{remark}
\newtheorem{remark}[theorem]{Remark}
\newcommand{\Hh}{\mathcal H}

\DeclareMathOperator{\CN}{CN}
\numberwithin{equation}{section}
\allowdisplaybreaks[1]
\date{}

\begin{document}
	
\title{Extremal 3-uniform hypergraphs without generalized 4-cycles
\thanks{ Leilei Zhang's work was supported by  the National Natural Science Foundation of China (Grant Nos. 12601664, 12571365), JSPS KAKENHI (Grant No.~25KF0036), the Hubei Provincial Natural Science Foundation of China (Grant No.~2025AFB309), the China Postdoctoral Science Foundation (Grant No.~2025M773113), the Fundamental Research Funds for the Central Universities, Central China Normal University (Grant No.~CCNU24XJ026). }}
	
\author[1, 2]{Leilei Zhang\footnote{Email: \texttt{mathdzhang@163.com}}}

\affil[1]{\footnotesize School of Mathematics and Statistics, Key Laboratory of Nonlinear Analysis \& Applications (Ministry of Education), Central China Normal University, Wuhan 430079, China}
\affil[2]{\small Faculty of Environment and Information Sciences, Yokohama National University, Yokohama 240-8501, Japan}

\maketitle

\vspace{-1.5em} 
	
\begin{abstract}

Let $f_r(n)$ be the maximum number of edges in an $n$-vertex $r$-uniform hypergraph in which all disjoint pairs of edges have distinct unions. In 1977, Erd\H{o}s posed the problem of determining $f_r(n)$ and, in particular, asked whether $f_3(n)=O(n^2)$. F\"uredi later proved that
$\binom{n-1}{r-1}+\left\lfloor (n-1)/r\right\rfloor
\le f_r(n)<\frac72\binom{n}{r-1}$
for $r\ge3$ and $n\ge2r$. Motivated by his constructions, he further conjectured that $f_3(n)\le\binom{n}{2}$ for all sufficiently large $n$, and that, for each fixed $r\ge4$, the lower bound is exact for all sufficiently large $n$. This is known as the Erd\H{o}s--F\"uredi conjecture. In a recent breakthrough, Huang, Ma and Yang proved the conjecture for every fixed $r\ge4$ and characterized all extremal hypergraphs. In this paper, we settle the remaining $3$-uniform case for every $n$. The bound is sharp whenever $n\ge5$ and $n\equiv1,5\pmod{20}$. Our proof uses a local weight inequality and a transfer argument.
	\end{abstract}
	
	\section{Introduction}
	
 In this paper, we complete the resolution of the Erd\H{o}s--F\"uredi conjecture on generalized $4$-cycles by settling its remaining 3-uniform case for every $n$. This conjecture arose from a hypergraph Tur\'an problem posed by Erd\H{o}s nearly half a century ago. A hypergraph $H$ consists of a vertex set $V(H)$ and an edge set $E(H)$. For an integer $r\ge2$, we call $H$ \emph{$r$-uniform}, or an \emph{$r$-graph}, if $E(H)\subseteq\binom{V(H)}r$. A \emph{generalized $4$-cycle} is an $r$-uniform hypergraph with four distinct edges $A,B,C,D$ satisfying $ A\cup B=C\cup D $ and $A\cap B=C\cap D=\varnothing. $ Let $\mathcal C_4^r$ denote the family of all such hypergraphs, and let $f_r(n)$ be the maximum number of edges in an $n$-vertex $r$-uniform hypergraph containing no member of $\mathcal C_4^r$. 
	
 When $r=2$, the family $\mathcal C_4^2$ consists, up to isomorphism,  of the ordinary $4$-cycle. Thus $f_2(n)=\operatorname{ex}(n,C_4)$, and determining this Tur\'an number is one of the classical
 problems in extremal graph theory. This problem dates back to Erd\H{o}s~\cite{Erdos3} in 1938. The classical results of Erd\H{o}s, R\'enyi and S\'os~\cite{ERS}, and independently Brown~\cite{brown}, yield $ f_2(n)=\left(\frac12+o(1)\right)n^{3/2}. $ F\"uredi~\cite{FurediC4} subsequently determined the exact value for an infinite sequence of orders, proving that $ f_2(q^2+q+1)=\frac12q(q+1)^2$  whenever $q=2^\alpha,\ \alpha\ge1.$ For more recent refinements of the upper bounds for $f_2(n)$, see Ma and Yang~\cite{MY}.
	
The corresponding problem for hypergraphs was posed by Erd\H{o}s~\cite{Erdos} in 1977, who asked whether $f_r(n)=O(n^{r-1})$ for every fixed $r\ge3$. He also reported an unpublished proof, jointly with Bollob\'as, that $f_3(n)=O(n^2)$. F\"uredi~\cite{Furedi} answered the general question by proving that, for $r\ge3$ and $n\ge2r$,
	\[
	\binom{n-1}{r-1}+\left\lfloor\frac{n-1}{r}\right\rfloor
	\le f_r(n)<\frac72\binom n{r-1}.
	\]
	
The lower bound follows by taking all $r$-sets containing a fixed vertex together with a maximum matching on the remaining vertices. For $r=3$, F\"uredi~\cite{Furedi} gave a stronger construction: replacing each block of a Steiner system $S(2,5,n)$ by all its $3$-element subsets yields a $\mathcal C_4^3$-free hypergraph with $\binom n2$ edges. By Hanani's theorem~\cite{Hanani}, such systems exist whenever $n\ge5$ and $n\equiv1,5\pmod{20}$, so $f_3(n)\ge\binom n2$ for these orders. These constructions led F\"uredi to propose the following conjecture~\cite{Furedi}, which was later included in a problem list of Erd\H{o}s~\cite{Erdos2}.
	
\begin{conjecture}[Erd\H{o}s--F\"uredi]\label{conj:EF}
For all sufficiently large $n$, $f_3(n)\le\binom n2$. Moreover, for every fixed $r\ge4$ and all sufficiently large $n$, $f_r(n)=\binom{n-1}{r-1}+\left\lfloor\frac{n-1}{r}\right\rfloor.$
\end{conjecture}
	
F\"uredi~\cite{Furedi} also showed that $f_3(n)/\binom n2$ converges as $n\to\infty$. Subsequent work progressively reduced the gap between the upper and lower bounds. Mubayi and Verstra\"ete~\cite{MVbipartite} improved the general upper bound to $ f_r(n)\le3\binom n{r-1}+O_r(n^{r-2}), (r\ge3). $ Pikhurko and Verstra\"ete~\cite{PV} further proved that
$$
	\limsup_{n\to\infty}\frac{f_r(n)}{\binom n{r-1}}
	\le\min\left\{\frac74,1+\frac2{\sqrt r}\right\}
	\qquad(r\ge3),
$$
and established the stronger estimate $f_3(n)\le\frac{13}{9}\binom n2$ for every $n\ge1$. Their approach extends F\"uredi's graph deletion lemma from bipartite graphs to general graphs and combines it with a randomized argument. Recently, Huang, Ma and Yang~\cite{HMY} resolved Conjecture~\ref{conj:EF} for every fixed $r\ge4$ and all sufficiently large $n$. They also determined all extremal hypergraphs, as stated below. We write $[n]=\{1,\ldots,n\}$.
	
\begin{theorem}[Huang--Ma--Yang~\cite{HMY}]
For every fixed $r\ge 4$ and all sufficiently large $n$,
	\[
	f_r(n)=\binom{n-1}{r-1}
	+\left\lfloor\frac{n-1}{r}\right\rfloor.
	\]
Moreover, every extremal $r$-graph is isomorphic to one of the following two constructions:
\begin{enumerate}
\item[(C1).] The $r$-graph with vertex set $[n]$ and edge set
 \[ 
 \left\{ A\in\binom{[n]}{r}:n\in A \right\}
		\cup
 \left\{ \{kr+1,kr+2,\ldots,kr+r\}: 0\le k\le\left\lfloor\frac{n-1}{r}\right\rfloor-1 \right\}.
 \]
		
\item[(C2).] When $r\mid n$, the $r$-graph obtained from (C1) by replacing the edge $\{n-r+1,\ldots,n\}$ with $\{n-r,\ldots,n-1\}$.
\end{enumerate}
\end{theorem}
	
Several partial results are known for the remaining case $r=3$. Up to isomorphism, $\mathcal C_4^3$ consists of a single hypergraph with edge set
\[
\bigl\{\{x,a,b\},\{y,a,b\},\{x,c,d\},\{y,c,d\}\bigr\},
\]
where $x,y,a,b,c,d$ are distinct vertices.

Mubayi~\cite{MubayiExact} proved that, for every $n$, an $n$-vertex $3$-uniform hypergraph with more than $\binom n2$ edges contains either a generalized $4$-cycle or a complete tripartite $3$-uniform hypergraph with part sizes $1,2,4$. Marchant~\cite{Marchant} proved that every $\mathcal C_4^3$-free $3$-uniform hypergraph with at least two vertices has a pair of vertices contained in at most three edges. Both bounds are sharp whenever $n\ge5$ and $n\equiv1,5\pmod{20}$: the design construction above has exactly $\binom n2$ edges, and every pair of vertices belongs to exactly three edges. Another related result concerns \emph{minimal $4$-cycles}, whose four edges admit a cyclic ordering in which two distinct edges intersect if and only if they are consecutive. Mubayi and Verstra\"ete~\cite{MVminimal} showed that, for every fixed $r\ge3$, an $n$-vertex $r$-uniform hypergraph containing no minimal $4$-cycle has at most $\binom{n-1}{r-1}+O_r(n^{r-2})$ edges. Every generalized $4$-cycle is a minimal $4$-cycle, so forbidding all minimal $4$-cycles imposes a stronger restriction.
	
In this paper, we settle the $3$-uniform hypergraph case of Conjecture~\ref{conj:EF} for every $n$. More precisely, we prove a stronger inequality involving two families of vertex pairs. For a $3$-uniform hypergraph $\Hh$ with vertex set $V$, let $\partial\Hh$ denote its shadow, consisting of all pairs contained in an edge of $\Hh$, and define
\[
\begin{aligned}
D(\Hh)=\Bigl\{\{x,y\}\in\binom V2:{} \{x,a,b\},\{y,a,b\}\in E(\Hh)\,\text{for some }\{a,b\}\in\binom{V\setminus\{x,y\}}2\Bigr\}.
\end{aligned}
\]
	
\begin{theorem}\label{thm:main}
Every $\mathcal C_4^3$-free $3$-uniform hypergraph $\Hh$ satisfies
\[2|E(\Hh)|\le|\partial\Hh|+|D(\Hh)|.\]	
\end{theorem}
	
Both $\partial\Hh$ and $D(\Hh)$ are contained in $\binom V2$, so Theorem~\ref{thm:main} immediately gives the following result.
	
\begin{theorem}\label{cor:main}
For every integer $n\ge1$, we have $f_3(n)\le\binom n2$. Equality holds whenever $n\ge5$ and $n\equiv1,5\pmod{20}$.
\end{theorem}

The design construction also attains equality in Theorem~\ref{thm:main}. Indeed, let $\Hh$ be obtained from a Steiner system $S(2,5,n)$ on $V$. Every pair $\{x,y\}\subseteq V$ lies in a unique block $B$ of size five and belongs to $\partial\Hh$. Choosing two distinct vertices $a,b\in B\setminus\{x,y\}$ gives $\{x,a,b\},\{y,a,b\}\in E(\Hh)$, so $\{x,y\}\in D(\Hh)$ as well. Thus $\partial\Hh=D(\Hh)=\binom V2$, and $|E(\Hh)|=\binom n2$ gives equality.
	
Our proof assigns a charge to each edge so that the total charge equals $|\partial\Hh|+|D(\Hh)|$, and then shows that the average charge is at least two. In Section~\ref{sec:charges}, we express these charges through auxiliary tripartite graphs. Section~\ref{sec:local-estimates} identifies the local structures whose charge is below two. Finally, Section~\ref{sec:proof} redistributes charge from edges with sufficient surplus to compensate for these deficits and completes the proof.

\section{Local charges}\label{sec:charges}
Throughout the proof, let $\Hh$ be a finite simple $3$-uniform hypergraph containing no four edges $\{x,a,b\},\{y,a,b\},\{x,c,d\},\{y,c,d\}$ on six distinct vertices. For a pair $p\in\binom{V(\Hh)}2$, let $N(p)=\{v\in V(\Hh)\setminus p:p\cup\{v\}\in E(\Hh)\}$ and $d_{\Hh}(p)=|N(p)|$. For distinct vertices $x,y$, their \emph{common link} $\CN(x,y)$ is the graph with edge set
\[ 
E(\CN(x,y))=\left\{\{u,v\}\in\binom{V(\Hh)\setminus\{x,y\}}2: \{x,u,v\},\{y,u,v\}\in E(\Hh)\right\}.
\]
Its vertex set is the union of these edges, so it has no isolated vertices. We also write $\CN(q)$ when $q=\{x,y\}$. The forbidden configuration is equivalent to two disjoint edges in a common link. Thus every common link is intersecting. A nonempty intersecting simple graph without isolated vertices is a star or a triangle: if it contains $\{a,b\},\{a,c\}$ and an edge avoiding $a$, that edge is $\{b,c\}$, and every edge meeting all three belongs to this triangle. A single edge is regarded as a star.
	
An edge $E\in E(\Hh)$ is a \emph{wing} of a \emph{root} $q=\{x,y\}$ if $E=\{x,u,v\}$ or $E=\{y,u,v\}$ for some $\{u,v\}\in E(\CN(q))$. Define its charge by
\[
	C(E)=\sum_{p\in\binom E2}\frac1{d_{\Hh}(p)} +\sum_{q:\,E\text{ is a wing of }q}\frac1{2|E(\CN(q))|}.
\]
Each pair in $\partial\Hh$ contributes one unit in total to its incident edges. Each root in $D(\Hh)$ contributes one unit to its $2|E(\CN(q))|$ distinct wings. Hence
\begin{equation}\label{eq:charge-sum}
	\sum_{E\in E(\Hh)}C(E)=|\partial\Hh|+|D(\Hh)|.
\end{equation}
	
For $E=\{a,b,c\}$, let
\[
\begin{aligned}
		X&=\{x_X:x\in N(\{a,b\})\setminus\{c\}\},\\
		Y&=\{y_Y:y\in N(\{a,c\})\setminus\{b\}\},\\
		Z&=\{z_Z:z\in N(\{b,c\})\setminus\{a\}\}.
\end{aligned}
\]
The subscripts distinguish copies of the same original vertex in different parts; all original labels lie outside $E$. Define a tripartite graph $K_E$ on $X,Y,Z$ by
\[
\begin{aligned}
		\{x_X,y_Y\}\in E(K_E)&\quad\Longleftrightarrow\quad \{a,x,y\}\in E(\Hh),\\
		\{x_X,z_Z\}\in E(K_E)&\quad\Longleftrightarrow\quad \{b,x,z\}\in E(\Hh),\\
		\{y_Y,z_Z\}\in E(K_E)&\quad\Longleftrightarrow\quad \{c,y,z\}\in E(\Hh).
\end{aligned}
\]
Vertices with equal labels are nonadjacent. For any tripartite graph $K$ with specified parts $V_1,V_2,V_3$, possibly empty, let
\[
\Phi(K)=\sum_{i=1}^3\frac1{|V_i|+1}+\frac12\sum_{v\in V(K)}\frac1{\deg_K(v)+1}.
\]
We call the first three terms the part contributions. A graph obtained by the construction of $K_E$ from $\Hh$, or from a subhypergraph containing $E$, is an \emph{exact tagged graph}. We will also evaluate $\Phi$ after adding auxiliary edges; such a graph need not arise from a $3$-uniform hypergraph. When the parts are denoted by $V_1,V_2,V_3$, the copy of label $x$ in $V_i$ is denoted by $x_i$.
	
For $x_X\in X$, the common link $\CN(c,x)$ contains $\{a,b\}$. Every other edge meets $\{a,b\}$ and has the form $\{a,y\}$ or $\{b,z\}$, corresponding to a neighbor $y_Y$ or $z_Z$ of $x_X$. Therefore $|E(\CN(c,x))|=1+\deg_{K_E}(x_X)$. The analogous roots represented by $Y,Z$ have $b,a$, respectively, as their unique vertex in $E$. These three classes are disjoint and include every root for which $E$ is a wing. Since the three pair degrees of $E$ are $|X|+1,|Y|+1,|Z|+1$, we obtain
\begin{equation}\label{eq:local-charge}
		C(E)=\Phi(K_E).
\end{equation}
	
A nonisolated tagged vertex is \emph{ordinary} if all its neighbors lie in one other part, and is a \emph{switch} otherwise. If $x_X$ has neighbors $y_Y,z_Z$, then $\CN(c,x)$ contains $\{a,b\},\{a,y\},\{b,z\}$. Intersection forces $y=z$ and excludes every further edge. Thus every switch has exactly two neighbors, one in each other part, with the same original label. We call these neighbors its \emph{targets} and their common label its \emph{target label}. Its common link is a triangle. Conversely, every triangular common link for which $E$ is a wing determines one switch. This correspondence is a bijection.
	
The tagged vertices $x_X,y_Y,z_Z$ represent the \emph{completion edges} $\{a,b,x\},\{a,c,y\},\{b,c,z\}$, respectively. Deleting any selection of completion edges deletes precisely the corresponding tagged vertices and leaves all edges between surviving vertices unchanged. Indeed, each completion edge contains two vertices of $E$, whereas an auxiliary edge is witnessed by a triple containing exactly one. Thus every induced vertex deletion is exact for a subhypergraph retaining $E$.
	
\begin{lemma}\label{lem:isolates-leaves}
Adding an isolated vertex to a part of size $h\ge0$ changes $\Phi$ by
\[
\frac12-\frac1{(h+1)(h+2)}\ge0.
\]
Adding a leaf to a part of size $h\ge1$, with its neighbor of degree $d\ge1$ before the addition, changes $\Phi$ by
\[
\frac14-\frac1{(h+1)(h+2)}-\frac1{2(d+1)(d+2)}\ge0.
\]
\end{lemma}
	
\begin{proof}
		Both formulas follow by subtraction. In the second formula, the two losses are at most $1/6$ and $1/12$.
	\end{proof}
	
	For positive integers $a,b$, let
	\[
	P(a,b)=\frac12\left(\frac a{b+1}+\frac b{a+1}\right),
	\qquad
	g(a,b)=\frac{(a-b)(a-b+1)}{2(a+1)(b+1)}.
	\]
	The first quantity is the vertex contribution of a complete bipartite graph. Direct calculation gives
	\begin{equation}\label{eq:P-g}
		P(a,b)+\frac1{b+1}=1+g(a,b),\qquad g(a,b)\ge0,
	\end{equation}
	and
	\begin{equation}\label{eq:P-balance}
		P(a,b)-\frac{a+b}{a+b+2}
		=\frac{(a+b+1)(a-b)^2}{2(a+1)(b+1)(a+b+2)}\ge0.
	\end{equation}
	The inequality $g(a,b)\ge0$ uses the fact that $a-b$ is an integer. In particular, $P(a,b)\ge1/2$.
	
	\begin{lemma}\label{lem:no-switch}
		Let $K$ be a tripartite graph in which every vertex has neighbors in at most one other part. If $M$ is its largest part size, then $\Phi(K)\ge2+1/(M+1)$. Consequently, if $E\in E(\Hh)$ is a wing of no triangular common link, then $C(E)>2$.
	\end{lemma}
	
	\begin{proof}
		First suppose that $K$ has no isolated vertices. Partition each part according to the part containing its vertices' neighbors, and add all edges between opposite classes. This cannot increase $\Phi$ and leaves at most three complete bipartite components. If no component remains, then $M=0$ and $\Phi=3$. For one component, with part sizes $a,c$,
		\[
		\Phi-2=\frac{(a-c)^2+a+c+2}{2(a+1)(c+1)}
		\ge\frac12\left(\frac1{a+1}+\frac1{c+1}\right)
		\ge\frac1{M+1}.
		\]
		For two components, with sizes $a,c$ and $b,e$, where the classes of sizes $a,b$ lie in the common part,~\eqref{eq:P-g} gives
		\[
		\Phi=2+\frac1{a+b+1}+g(a,c)+g(b,e)\ge2+\frac1{M+1}.
		\]
		
		For three components, remove one with sizes $d,f$. Let $c,e\ge1$ be the remaining sizes in those two parts. Restoring the removed component changes $\Phi$ by
		\[
		\Delta=P(d,f)-\frac d{(c+1)(c+d+1)}-\frac f{(e+1)(e+f+1)}.
		\]
		Since $c,e\ge1$,
		\[
		2\Delta\ge\frac d{f+1}+\frac f{d+1}-\frac d{d+2}-\frac f{f+2}>0.
		\]
		Indeed, replacing the last two denominators by $d+1,f+1$ makes the right-hand side $(d-f)^2/((d+1)(f+1))\ge0$. The common part of the two retained components is unchanged during restoration and has size at most $M$. Hence their lower bound proves the result for three components as well.
		
		For a graph with isolated vertices, first remove them and apply the preceding argument. Restoring them cannot decrease $\Phi$, by Lemma~\ref{lem:isolates-leaves}, and can only increase the largest part size. The last assertion follows from~\eqref{eq:local-charge} and the correspondence between switches and triangular common links.
	\end{proof}
	
	\begin{lemma}\label{lem:labels}
		If $x_X$ is a switch with targets $y_Y,y_Z$, then
		\[
		\{z:z_Y\in Y,\ z_Z\in Z\}\subseteq\{x,y\}.
		\]
		If distinct labels $r,s$ occur in both $V_i,V_j$, then $\{r_i,s_j\}\in E(K_E)$ if and only if $\{s_i,r_j\}\in E(K_E)$. Two switches in one part with different target labels $u,v$ have source labels $v,u$, respectively. A part with at least three switches has one common target label, and both targets are ordinary. More generally, two or more switches in one part sharing one target share both targets, and these targets are ordinary.
	\end{lemma}
	
	\begin{proof}
		The switch gives $\{a,x,y\},\{b,x,y\}\in E(\Hh)$, so $\{x,y\}\in E(\CN(a,b))$. If $z$ occurs in $Y,Z$, then $\{a,c,z\},\{b,c,z\}\in E(\Hh)$, so $\{c,z\}\in E(\CN(a,b))$. Intersection forces $z\in\{x,y\}$. The two edges in the second assertion are witnessed by the same triple.
		
		Suppose two switches have labels $r,s$ and distinct target labels $u,v$. The first assertion gives $\{u,v\}\subseteq\{r,u\}\cap\{s,v\}$, so $r=v$ and $s=u$. If the other two parts share two labels, the same restriction allows at most two switch labels in the given part. Thus three switches have one common target label. Each target then has at least three neighbors in that part and cannot be a switch. Finally, switches sharing one target have the same target label, so they also share the other target. Each target has at least two neighbors in the same part and is therefore ordinary.
	\end{proof}
	
	Let $t_i$ be the number of switches in $V_i$. An ordinary vertex belongs to the class $i\to j$ if its neighbors lie in $V_j$; denote the size of this class by $a_{ij}$. We call $\{i,j\}$ \emph{active} if $a_{ij},a_{ji}>0$. Adding all edges between opposite ordinary classes only decreases $\Phi$. We use this completion solely for numerical estimates; all arguments involving labels refer to the exact graph before completion.
	
	\begin{lemma}\label{lem:pruning}
		An exact tagged graph can be reduced by induced deletions to one with at most two switches in each part, without increasing $\Phi$, creating isolated vertices, or changing its ordinary classes and active pairs. Each part originally containing at least two switches retains at least two.
	\end{lemma}
	
	\begin{proof}
		If a part contains $m\ge3$ switches, Lemma~\ref{lem:labels} gives two common ordinary targets. Let $h$ be its size and $d_1,d_2$ the target degrees. All three numbers are at least $m$. Deleting one switch $s$ gives
		\[
		\begin{aligned}
			\Phi(K)-\Phi(K-s)
			&=\frac16-\frac1{h(h+1)}-\frac1{2d_1(d_1+1)}-\frac1{2d_2(d_2+1)}\\
			&\ge\frac16-\frac2{m(m+1)}\ge0.
		\end{aligned}
		\]
		Both targets remain ordinary with at least two neighbors in that part. No other vertex changes degree. Since only a switch is deleted and the targets remain in their ordinary classes, all ordinary classes and active pairs are unchanged. Repeat until the part has two switches, and apply the same argument to any other part with at least three. The deletions are realized by deleting completion edges, so every remaining graph is exact.
	\end{proof}
	
	\begin{lemma}\label{lem:six-cycle}
		If a part contains at least two switches and one has a switch neighbor, then all switches form a six-cycle, with two in each part. Otherwise, every part incident with a switch--switch edge has exactly one switch.
	\end{lemma}
	
	\begin{proof}
		Suppose $\{r_1,s_2\}$ joins switches. Their other targets are $s_3,r_3$. A second switch in $V_1$ cannot have target label $s$, since $s_2$ would then have two neighbors in $V_1$. By Lemma~\ref{lem:labels}, this second switch is $s_1$, with targets $r_2,r_3$. Applying the same lemma to the edge $\{s_2,r_3\}$ gives the edge $\{r_2,s_3\}$. Thus the six vertices are switches on a cycle in the order $r_1,s_2,r_3,s_1,r_2,s_3,r_1$. Each part has two different target labels, so Lemma~\ref{lem:labels} excludes every further switch.
	\end{proof}
	
	\begin{lemma}\label{lem:inactive}
		Let $h_{ij}$ be the number of switch--switch edges between $V_i,V_j$. These edges form a matching. Switches send $t_j-h_{ij}$ edges to the ordinary class $i\to j$, and $t_i-h_{ij}$ edges to $j\to i$. If $\{i,j\}$ is inactive, then $h_{ij}=\min\{t_i,t_j\}$. If also $t_i=t_j$, both ordinary classes are empty. Otherwise exactly one is nonempty, in the part with fewer switches; its size $r$ satisfies $1\le r\le|t_i-t_j|$, and its vertex contribution is at least
		\[
		\frac{r^2}{2(r+|t_i-t_j|)}.
		\]
	\end{lemma}
	
	\begin{proof}
		Every switch has one neighbor in each other part, which proves the matching assertion and the edge counts. For an inactive pair, one ordinary class is empty, so the matching saturates the smaller switch class. The remaining $|t_i-t_j|$ edges meet precisely the ordinary class in the smaller part. Its vertices have positive degree and total degree $|t_i-t_j|$. Cauchy's inequality gives
		\[
		\frac12\sum_{v\text{ in the class}}\frac1{\deg_K(v)+1}
		\ge\frac{r^2}{2(r+|t_i-t_j|)}.
		\]
	\end{proof}
	
	Each incidence of a switch is called a \emph{port} and receives $1/12$ of the switch's contribution $1/6$. Adding a switch neighbor to a completed active ordinary class reduces a vertex contribution by
	\[
	\frac1{2(d+1)}-\frac1{2(d+2)}
	=\frac1{2(d+1)(d+2)}\le\frac1{12}\qquad(d\ge1).
	\]
	The corresponding port covers this loss. An ordinary vertex of degree $s\ge1$ on an inactive pair, together with its incident ports, contributes at least $1/3$, since
	\[
	\frac s{12}+\frac1{2(s+1)}-\frac13
	=\frac{(s-1)(s-2)}{12(s+1)}\ge0.
	\]
	Hence, if $R$ is the number of ordinary vertices on inactive pairs, then
	\begin{equation}\label{eq:ports}
		\Phi(K)\ge\sum_{i=1}^3\frac1{|V_i|+1}
		+\sum_{\substack{i<j\\\{i,j\}\text{ active}}}P(a_{ij},a_{ji})+\frac R3.
	\end{equation}
	Every ordinary vertex belongs to one class, and each port is used at most once. The remaining contributions, from ports on switch--switch edges and from isolated vertices, are nonnegative.
	
	\section{Local estimates}\label{sec:local-estimates}
	
	We first record the numerical inequalities used to estimate the losses in the part contributions. Recall the functions $P$ and $g$ from Section~2.
	
	\begin{lemma}\label{lem:numerical}
		For positive integers $a,b,c$, let $L(c)=1/((c+1)(c+2))$ and $B(c)=3/((c+1)(c+3))$. Then
		\begin{equation}\label{eq:numerical-loss}
			L(c)\le g(a,c)+\frac1{2(a+2)},\qquad
			B(c)\le g(a,c)+\frac3{2(a+2)},
		\end{equation}
		and
		\begin{equation}\label{eq:numerical-sum}
			\begin{aligned}
				\frac1{2(a+2)}+\frac1{2(b+2)}&\le\frac16+\frac1{a+b+3},\\
				\frac3{2(a+2)}+\frac1{2(b+2)}&\le\frac12+\frac1{a+b+3}.
			\end{aligned}
		\end{equation}
	\end{lemma}
	
	\begin{proof}
		If $a\notin\{c,c-1\}$, then $g(a,c)\ge L(c)$. Indeed, if $a=c+s$ with $s\ge1$, use $a+1\le s(c+2)$ in the definition of $g$. If $a\le c-2$, use $(a-c)(a-c+1)\ge2$ and $a+1<c+2$. In the remaining two cases, $a\le c$ and $1/(2(a+2))\ge L(c)$. This proves the first inequality in~\eqref{eq:numerical-loss}.
		
		For the second, if $a\le c+1$, then $3/(2(a+2))\ge B(c)$. Otherwise $a=c+s$ for some $s\ge2$. The inequalities $s(s+1)\ge3s$ and $s/(c+s+1)\ge2/(c+3)$ give $g(a,c)\ge B(c)$.
		
		For the first inequality in~\eqref{eq:numerical-sum}, assume $a\le b$. The two summands are at most $1/6$ and $1/(a+b+3)$, respectively. For the second, the corresponding bounds $1/2$ and $1/(a+b+3)$ hold when $a\le b+1$. Otherwise $a\ge b+2\ge3$, and the sum is at most $3/10+1/6=7/15<1/2$.
	\end{proof}
	
	\begin{lemma}\label{lem:two-active}
		If an exact tagged graph $K$ has no isolated vertices and at least two active ordinary pairs, then $\Phi(K)\ge2$.
	\end{lemma}
	
	\begin{proof}
		By Lemma~\ref{lem:pruning}, we may assume that each part contains at most two switches. Let $T\le6$ be their total number.
		
		Suppose first that all three pairs are active. Let $u_1,u_2,u_3\ge2$ be the total ordinary class sizes on the three pairs, and let $U=u_1+u_2+u_3$. By~\eqref{eq:P-balance},
		\[
		\sum_{i<j}P(a_{ij},a_{ji})
		\ge\sum_{i=1}^3\frac{u_i}{u_i+2}
		\ge2-\frac2{U-2}.
		\]
		For the last inequality, the concavity of $u/(u+2)$ shows that, for fixed $U$, the sum is minimized when two variables equal two. The part sizes sum to $U+T$. Hence~\eqref{eq:ports} and Cauchy's inequality yield
		\[
		\Phi(K)\ge2-\frac2{U-2}+\frac9{U+T+3}\ge2,
		\]
		where the last inequality is equivalent to $7U\ge2T+24$ and follows from $U\ge6$ and $T\le6$.
		
		Now suppose that exactly two pairs are active. Denote them by $\{0,1\}$ and $\{0,2\}$, and let their ordinary class sizes be $a,c$ and $b,e$, respectively, with the $a$- and $b$-classes in part $0$. Let $t_0,t_1,t_2\le2$ be the switch counts in the three parts. If the switches form a six-cycle, none meets an ordinary vertex. Their contribution is one, and each ordinary pair contributes at least $1/2$ by~\eqref{eq:P-balance}. Thus $\Phi(K)\ge2$.
		
		Assume that the switches do not form a six-cycle. Since $\{1,2\}$ is inactive, Lemma~\ref{lem:inactive} gives $h_{12}=\min\{t_1,t_2\}$. If both counts are positive and one exceeds one, Lemma~\ref{lem:six-cycle} gives the excluded six-cycle. Therefore, up to interchanging parts $1$ and $2$,
		\[
		(t_1,t_2)\in\{(0,0),(1,0),(2,0),(1,1)\}.
		\]
		Let $r_i$ be the size of the ordinary class on the inactive pair in part $i$. If $t_1=t_2$, then $r_1=r_2=0$. If $(t_1,t_2)=(q,0)$, then $r_1=0$ and $1\le r_2\le q$.
		
		At most $2t_0+t_1+t_2$ switch ports meet active ordinary classes. Paying the loss at each such attachment with $1/12$ leaves switch contribution at least $(t_1+t_2)/12$. For $c\ge1$ and $k\ge0$,
		\[
		\frac1{c+1}-\frac1{c+k+1}\le\frac12-\frac1{k+2}.
		\]
		Using this inequality for the outer part contributions,~\eqref{eq:P-g} for the two active pairs, and Lemma~\ref{lem:inactive} for the inactive classes, we obtain
		\[
		\Phi(K)\ge1+\frac1{t_1+r_1+2}+\frac1{t_2+r_2+2}
		+\frac{t_1+t_2}{12}
		+\sum_{i=1}^2\frac{r_i^2}{2(r_i+|t_1-t_2|)}.
		\]
		A summand with zero numerator is interpreted as zero. We have discarded the central part contribution and the nonnegative $g$ terms. For $(t_1,t_2)=(1,0)$, necessarily $r_2=1$, and the terms after the initial one sum to one. For $(t_1,t_2)=(2,0)$ and $r_2=2$, they sum to $7/6$. Both cases give $\Phi(K)\ge2$.
		
		If $(t_1,t_2)=(0,0)$, retaining the discarded terms gives
		\[
		\Phi(K)\ge2+\frac1{a+b+t_0+1}+g(a,c)+g(b,e)>2.
		\]
		If $(t_1,t_2)=(1,1)$, the two outer part contributions lose $L(c)$ and $L(e)$. Thus
		\[
		\Phi(K)\ge2+\frac1{a+b+t_0+1}+g(a,c)+g(b,e)
		+\frac16-L(c)-L(e)\ge2.
		\]
		The last inequality follows from Lemma~\ref{lem:numerical} and $t_0\le2$.
		
		It remains to consider $(t_1,t_2)=(2,0)$ and $r_2=1$. The two switches meet the same inactive ordinary vertex. By Lemma~\ref{lem:labels}, their other targets also coincide. This target is ordinary, since a switch target would give a six-cycle by Lemma~\ref{lem:six-cycle}; it therefore belongs to the $a$-class. Its two attachments lose exactly $1/((c+1)(c+3))$ from its completed contribution. The part $1$ reciprocal loses $2/((c+1)(c+3))$. The inactive ordinary vertex contributes $1/6$, and its addition to part $2$ loses $L(e)$. The two outer switches contribute $1/3$, while the $t_0$ central switches pay for their own active attachments. Consequently,
		\[
		\Phi(K)\ge2+\frac1{a+b+t_0+1}+g(a,c)+g(b,e)
		+\frac12-B(c)-L(e)\ge2,
		\]
		again by Lemma~\ref{lem:numerical} and $t_0\le2$.
	\end{proof}
	
	\begin{lemma}\label{lem:no-active}
		If an exact tagged graph $K$ has no isolated vertices and no active ordinary pair, then $\Phi(K)\ge2$.
	\end{lemma}
	
	\begin{proof}
		Apply Lemma~\ref{lem:pruning}. Every ordinary vertex is a switch neighbor. If all switches lie in one part, there are one or two of them. For two switches, Lemma~\ref{lem:labels} shows that they either have the same target label or have two different, exchanged target labels. The resulting graphs are $K_{1,2}$, $K_{2,2}$ and $2K_{1,2}$, respectively. Their part sizes are $(1,1,1)$, $(2,1,1)$ and $(2,2,2)$.
		
		If switches lie in more than one part, Lemma~\ref{lem:inactive} forces a switch--switch edge between any two occupied parts. Lemma~\ref{lem:six-cycle} then gives a six-cycle unless every occupied part contains exactly one switch. Three such switches would be pairwise adjacent. The rule that both targets of a switch have the same original label would then give equal labels to adjacent vertices, which is impossible. Two switches and their ordinary targets form a four-vertex path, with part sizes $(1,1,2)$. Direct calculation gives
		\[
		\begin{array}{c|ccccc}
			K&K_{1,2}&K_{2,2}&2K_{1,2}&P_4&C_6\\ \hline
			\Phi(K)&13/6&2&7/3&13/6&2.
		\end{array}
		\]
		If there is no switch, there is no vertex and $\Phi(K)=3$. This proves the lemma.
	\end{proof}
	
	For the case of one active pair, we use two further inequalities. For $a,b\ge1$ and $0\le u,v\le2$,
	\begin{equation}\label{eq:one-pair-baseline}
		P(a,b)+\frac1{a+u+1}+\frac1{b+v+1}\ge1.
	\end{equation}
	Indeed,
	\[
	P(a,b)=1-\frac1{2(a+1)}-\frac1{2(b+1)}
	+\frac{(a-b)^2}{2(a+1)(b+1)},
	\]
	and $1/(a+3)\ge1/(2(a+1))$. Also, for positive integers $x,y$,
	\begin{equation}\label{eq:integer-identity}
		\frac{x+1}{y}+\frac{y-2}{x}
		=2-\frac1x+\frac{(y-x)(y-x-1)}{x\,y}
		\ge2-\frac1x.
	\end{equation}
	
	\begin{lemma}\label{lem:one-active}
		Let $K$ be an exact tagged graph with no isolated vertices and exactly one active ordinary pair. Then $\Phi(K)>23/12$. If $\Phi(K)<2$, then $K$ has exactly one switch, this switch lies in a part of the active pair, and the third part consists of one ordinary vertex. In this case, let $a,b\ge1$ be the ordinary class sizes, with the switch in the part containing the $a$-class. If $v$ is its target in the $b$-class and $w=1/(2(\deg_K(v)+1))$, then
		\begin{equation}\label{eq:target-charge}
			\Phi(K)-w\ge\frac{23}{12}+\frac{a}{2(a+1)(a+2)}.
		\end{equation}
	\end{lemma}
	
	\begin{proof}
		By Lemma~\ref{lem:pruning}, first assume $t_i\le2$. Let $\{1,2\}$ be the active pair, with ordinary class sizes $a,b\ge1$. Completing this pair will only be used to obtain numerical lower bounds. Up to interchanging parts $1$ and $2$, the cases and bounds proved below are as follows. In the seventh row, $m,n\in\{1,2\}$.
        \begin{center}
        \small
        \renewcommand{\arraystretch}{1.15}
        \begin{tabular}{p{44mm} p{56mm} c}
        \hline
        Switch counts & Structure of $V_3$ & Bound for $\Phi(K)$\\
        \hline
        $(0,0,0)$ & Empty & $>2$\\
        $(2,2,2)$, six-cycle & Two switches & $\ge2$\\
        $(0,0,1)$ & One switch & $\ge13/6$\\
        $(0,0,2)$ & Two switches & $\ge2$\\
        $(1,0,1)$ & One switch & $\ge2$\\
        $(1,1,1)$ & One switch & $\ge2$\\
        $(m,n,0)$ & Two ordinary vertices, of degrees $m,n$ & $\ge25/12$\\
        $(1,0,0)$ & One ordinary vertex, of degree one & $>23/12$\\
        $(2,0,0)$, same target label & One ordinary vertex, of degree two & $\ge2$\\
        $(2,0,0)$, different target labels & Two ordinary vertices, of degree one & $\ge13/6$\\
        \hline
        \end{tabular}
        \end{center}
        Only the eighth row can have charge below two; in that case we prove the stronger bound~\eqref{eq:target-charge}.

        If the switches form a six-cycle, their contribution is one and none meets an ordinary vertex. Inequality~\eqref{eq:one-pair-baseline} gives another one. Hence assume that this case does not occur.
		
		Suppose first that $t_3>0$. Both $\{1,3\}$ and $\{2,3\}$ are inactive, so Lemma~\ref{lem:inactive} gives $h_{i3}=\min\{t_i,t_3\}$ for $i=1,2$. If $t_3=2$, either $t_1=t_2=0$ or Lemma~\ref{lem:six-cycle} gives a six-cycle. If $t_3=1$, the same lemma forces $t_1,t_2\le1$. Thus, up to interchanging parts $1$ and $2$, the possible switch counts are
		\[
		(t_1,t_2,t_3)=(0,0,1),\ (0,0,2),\ (1,0,1),\ (1,1,1).
		\]
		
		For the first three cases, the two active part sizes exceed $a,b$ by at most two. By~\eqref{eq:one-pair-baseline}, their part contributions and the completed ordinary contribution are at least one, before accounting for switch attachments. For $(0,0,1)$, the two inactive ordinary vertices have degree one, the switch contributes $1/6$, and the third-part reciprocal is $1/2$. These give another $7/6$, with no active attachment.
		
		For $(0,0,2)$, each inactive ordinary class has size $r\in\{1,2\}$ and contributes at least $r^2/(2(r+2))\ge1/6$ by Lemma~\ref{lem:inactive}. The two switches and the third-part reciprocal contribute $2/3$. There is again no active attachment, so $\Phi(K)\ge2$.
		
		For $(1,0,1)$, the two switches are adjacent. Their other targets are ordinary vertices in part $2$, one in each of its two classes. The inactive target has degree one. This target, the two switches and the third-part reciprocal contribute $1/4+1/3+1/2=13/12$. There is one active attachment, whose loss is at most $1/12$. Thus $\Phi(K)\ge2$ in this case as well.
		
		For $(1,1,1)$, let the switch in part $3$ have label $y$. Its two switch neighbors have a common label $x$. Their other neighbors are the ordinary copies $y_1,y_2$, since adjacency between the two switches with label $x$ is impossible. Each active class therefore receives one attachment.
		
		More generally, suppose that groups of $m,n$ switches in the two active parts attach to one ordinary target each in the opposite active class. After completing the active pair, its ordinary contribution and the two part contributions equal
		\begin{equation}\label{eq:switch-groups}
			\begin{aligned}
				&\frac3{2(a+m+1)}+\frac3{2(b+n+1)}
				+\frac{a-1}{2(b+1)}+\frac{b-1}{2(a+1)}\\
				&\hspace{12mm}=1+\frac{(a-b)^2}{2(a+1)(b+1)}
				+f_m(a+1)+f_n(b+1),
			\end{aligned}
		\end{equation}
		where $f_t(x)=3/(2(x+t))-1/x$. Here $m=n=1$, and $f_1(x)=(x-2)/(2x(x+1))\ge0$ for $x\ge2$. The three switches and the third-part reciprocal contribute one more. Hence $\Phi(K)\ge2$. The completion may add the pair $\{y_1,y_2\}$ to the auxiliary edge set as a numerical relaxation; this pair need not belong to $E(K)$.
		
		Now suppose that $t_3=0$ and both active parts contain switches. We show that the switches in either part have a common target label. Suppose instead that part $1$ contains switches with different target labels. By Lemma~\ref{lem:labels}, their source labels are exchanged; denote these two labels by $r,s$. Their copies in part $2$ must be ordinary, since a switch neighbor would give a six-cycle by Lemma~\ref{lem:six-cycle}. Both $r$ and $s$ also occur in parts $1$ and $3$. Applying Lemma~\ref{lem:labels} to any switch in part $2$ forces its source label to be $r$ or $s$. This is impossible because both copies in that part are already ordinary. The argument for part $2$ is the same.
		
		Let the group sizes be $m,n\in\{1,2\}$. Their targets in part $3$ are two distinct ordinary vertices: one has neighbors only in part $1$, and the other only in part $2$. These are all the vertices in part $3$, by Lemma~\ref{lem:inactive}.

        We now justify applying~\eqref{eq:switch-groups} even when a group's target in the other active part is a switch. If that target is ordinary, completing the active pair leaves one vertex in the opposite ordinary class with the additional switch neighbors. If that target is a switch, the opposite ordinary class has no such attachments. In the latter case, choose any vertex of that ordinary class and replace its contribution by the smaller value it would have if those attachments were present. Consequently,
        \[
        \begin{aligned}
        \frac12\sum_{\substack{v\in V_1\cup V_2\\v\text{ ordinary}}}
             \frac1{\deg_K(v)+1}
        &\ge \frac{a-1}{2(b+1)}+\frac1{2(b+n+1)}\\
        &\quad+\frac{b-1}{2(a+1)}+\frac1{2(a+m+1)}.
        \end{aligned}
        \]
        Here completion can only decrease the left-hand contribution, and the additional replacement can only decrease it further. This is solely a numerical lower bound: we do not apply Lemmas~\ref{lem:labels} or~\ref{lem:inactive}, or any other structural assertion, to the virtual attachments. All switch contributions retain their actual value $1/6$. Adding the two active part contributions to the displayed lower bound gives the expression in~\eqref{eq:switch-groups}.
		
		For $x\ge2$,
		\[
		f_1(x)\ge0,\qquad
		f_2(x)+\frac18=\frac{(x-2)(x+8)}{8x(x+2)}\ge0.
		\]
		Thus $f_t(x)+t/6+1/(2(t+1))$ is at least $5/12$ for $t=1$ and at least $3/8$ for $t=2$. Including the third-part reciprocal $1/3$ in~\eqref{eq:switch-groups}, we obtain $\Phi(K)\ge4/3+2(3/8)=25/12$.
		
		Suppose next that only part $1$ contains switches, and let their number be $m\in\{1,2\}$. If they have a common target label, part $3$ consists of one ordinary vertex of degree $m$. The two active part contributions and their ordinary vertices are at least
		\[
		\frac3{2(a+m+1)}+\frac{a+2}{2(b+1)}+\frac{b-1}{2(a+1)}
		\ge1-\frac1{2(a+1)}+\frac3{2(a+m+1)},
		\]
		by~\eqref{eq:integer-identity} with $x=a+1$ and $y=b+1$. If $m=2$, the right-hand side is at least one. The third-part reciprocal, its ordinary vertex and the two switches contribute $1/2+1/6+2/6=1$, so $\Phi(K)\ge2$.
		
		If $m=1$, subtract the actual contribution $w$ of the target in the $b$-class before completing the active pair. All remaining summands can only decrease on completion. Therefore
		\[
		\begin{aligned}
			\Phi(K)-w
			&\ge\frac{11}{12}+\frac1{a+2}
			+\frac{a+2}{2(b+1)}+\frac{b-1}{2(a+1)}\\
			&\ge\frac{23}{12}+\frac{a}{2(a+1)(a+2)},
		\end{aligned}
		\]
		where the second inequality again follows from~\eqref{eq:integer-identity}. This proves~\eqref{eq:target-charge} and $\Phi(K)>23/12$ in this case.
		
		If $m=2$ and the target labels differ, part $3$ contains two ordinary vertices of degree one. The two targets in the $b$-class are distinct, so $b\ge2$. The switches, the third-part reciprocal and the two vertices in that part contribute $7/6$. For the remaining contribution, let $x=a+1$ and $y=b+1$. Completing the active pair and applying~\eqref{eq:integer-identity} gives the lower bound
		\[
		\begin{aligned}
			\frac{x+1}{2y}+\frac{y-3}{2x}+\frac1{x+1}+\frac1{x+2}
			&\ge1-\frac1x+\frac1{x+1}+\frac1{x+2}\\
			&\ge1,
		\end{aligned}
		\]
		since $-1/x+1/(x+1)+1/(x+2)=(x^2-2)/(x(x+1)(x+2))\ge0$ for $x\ge2$. Finally, if there is no switch, Lemma~\ref{lem:no-switch} applies.
		
		These cases prove the asserted bounds after pruning. If the original graph had charge below two, its pruned graph would also have charge below two. The only possible case has exactly one switch. However, every part from which switches were deleted retains two switches, by Lemma~\ref{lem:pruning}. Hence no pruning occurred. The structural description and~\eqref{eq:target-charge} therefore hold for the original graph.
	\end{proof}
	
	\begin{lemma}\label{lem:deficient}
		Every exact tagged graph satisfies $\Phi(K_E)>23/12$. If $C(E)<2$, then $K_E$ has no isolated vertex, exactly one active ordinary pair and exactly one switch. This switch lies in a part of the active pair, and the third part consists of one ordinary vertex. With $v,w$ as in Lemma~\ref{lem:one-active},
		\[
		\delta(E):=2-C(E)<\frac1{12}-w,\qquad \deg_{K_E}(v)\ge6.
		\]
		Moreover, $E$ contains exactly one pair of codegree two. In its unique triangular common link, this pair joins the vertex of the root lying in $E$ to a vertex of the triangle.
	\end{lemma}
	
	\begin{proof}
		Delete all isolated tagged vertices. As observed in Section~2, this is an exact operation, obtained by deleting the corresponding completion edges. Lemmas~\ref{lem:two-active}, \ref{lem:no-active} and~\ref{lem:one-active} give charge greater than $23/12$ for the remaining graph. Restoring isolated vertices cannot decrease the charge, by Lemma~\ref{lem:isolates-leaves}. This proves the universal bound.
		
		Suppose that $\Phi(K_E)<2$. By Lemma~\ref{lem:no-switch}, $K_E$ has a switch, so all three parts remain nonempty after deleting isolated vertices. Restoring an isolated vertex to a nonempty part increases the charge by at least $1/3$, by Lemma~\ref{lem:isolates-leaves}. Since $23/12+1/3>2$, there was no isolated vertex. Lemmas~\ref{lem:two-active}, \ref{lem:no-active} and~\ref{lem:one-active} now give the stated structure. With $a$ as in Lemma~\ref{lem:one-active}, inequality~\eqref{eq:target-charge} gives
		\[
		0<\delta(E)\le\frac1{12}-w-\frac{a}{2(a+1)(a+2)}
		<\frac1{12}-w.
		\]
		Thus $w<1/12$, and hence $\deg_{K_E}(v)\ge6$.
		
		The three part sizes are $a+1,b,1$. If $b=1$, the graph is obtained from a three-vertex $K_{1,2}$, with one vertex in each part, by adding leaves to the switch's part. Lemma~\ref{lem:isolates-leaves} gives $\Phi(K_E)\ge13/6$, a contradiction. Thus $b\ge2$, and exactly one pair of $E$ has codegree two.
		
		To identify its position, write $E=\{p,q,r\}$ and let the singleton part correspond to $N(\{p,q\})\setminus\{r\}=\{x\}$. Let the switch have label $z\in N(\{q,r\})\setminus\{p\}$. Its other target is the copy of $x$ in the part corresponding to $N(\{p,r\})\setminus\{q\}$. Hence $E(\CN(p,z))=\{\{q,r\},\{q,x\},\{r,x\}\}$. The codegree-two pair $\{p,q\}$ joins the root vertex $p$ to a vertex of this triangle, as required.
	\end{proof}
	
	\begin{remark}\label{rem:local-sharpness}
The constant $23/12$ in Lemma~\ref{lem:deficient} is the best possible infimum, and individual edges can have charge below two. To see this, take distinct vertices $a,b,c,r$ and disjoint sets $A,B$, also disjoint from these four vertices, with $|A|=|B|=m\ge1$. Fix $u\in B$. Let $\Hh_m$ have vertex set $\{a,b,c,r\}\cup A\cup B$ and edge set
\[
\begin{aligned}
E(\Hh_m)={}&\bigl\{\{a,b,c\},\{a,b,r\},\{b,c,u\},\{a,r,u\},\{b,r,u\}\bigr\}\\
&\cup\bigl\{\{a,b,x\}:x\in A\bigr\}
\cup\bigl\{\{a,c,y\}:y\in B\bigr\}\\
&\cup\bigl\{\{a,x,y\}:x\in A,\ y\in B\bigr\}.
\end{aligned}
\]
The only edges avoiding $a$ are $\{b,c,u\}$ and $\{b,r,u\}$. If two distinct disjoint pairs of edges had the same union, each pair would contain exactly one of these two edges. The edge paired with $\{b,c,u\}$ would then contain both $a$ and $r$. However, the only such edges, $\{a,b,r\}$ and $\{a,r,u\}$, both meet $\{b,c,u\}$. Thus $\Hh_m$ is $\mathcal C_4^3$-free.

For $E=\{a,b,c\}$, the three parts of $K_E$ are tagged copies of $A\cup\{r\}$, $B$ and $\{u\}$. The copies of $A$ and $B$ induce a complete bipartite graph, and the copy of $r$ is a switch whose targets are the copies of $u$ in the other two parts. The $m$ copies of $A$ have degree $m$, the $m-1$ copies of $B\setminus\{u\}$ have degree $m$, and the remaining copy in the second part has degree $m+1$. The switch has degree two and the third-part vertex has degree one. Hence
\[
\begin{aligned}
C(E)
&=\frac1{m+2}+\frac1{m+1}+\frac12+\frac16+\frac14
  +\frac{2m-1}{2(m+1)}+\frac1{2(m+2)}\\
&=\frac{23}{12}+\frac{2m+1}{2(m+1)(m+2)}
\longrightarrow\frac{23}{12}.
\end{aligned}
\]
In particular, $m=10$ gives $C(E)=527/264<2$, so the transfer argument is necessary.
\end{remark}

\begin{lemma}\label{lem:isolate-surplus}
		If all three parts of $K_E$ have size at least two and $K_E$ has $I$ isolated vertices, then $\Phi(K_E)\ge2+I/6$.
	\end{lemma}
	
	\begin{proof}
		If $I=0$, then Lemma~\ref{lem:deficient} gives $\Phi(K_E)\ge2$, since a graph of charge below two has a singleton part. Assume $I>0$, and delete all isolated vertices. The remaining graph is exact.
		
		If every remaining part is nonempty, restoring each isolated vertex adds at least $1/3$. Lemma~\ref{lem:deficient} therefore gives
		\[
		\Phi(K_E)>\frac{23}{12}+\frac I3\ge2+\frac I6.
		\]
		If some remaining part is empty, the remaining graph has no switch and has charge at least two by Lemma~\ref{lem:no-switch}. Restoring $k$ isolated vertices to an empty part adds
		\[
		\frac k2+\frac1{k+1}-1=\frac{k(k-1)}{2(k+1)}\ge\frac k6,
		\]
		since its original size was $k\ge2$. Restoring an isolated vertex to any other part adds at least $1/3$. Summing these contributions proves the result.
	\end{proof}
	
	\section{Proof of the main theorem}\label{sec:proof}

If $E(\CN(u,v))=\binom{T}{2}$ for a three-element set $T$, call $u,v$ its \emph{apices} and $T$ its \emph{base}. Its three wings at $u$ contain the three pairs of $T$. The pairs joining $u$ to the vertices of $T$ are its \emph{spokes at $u$}. Each spoke has codegree at least two; call it \emph{low} if its codegree is two, and \emph{high} otherwise. An edge $E\in E(\Hh)$ is \emph{deficient} if $C(E)<2$, and its deficit is $\delta(E)=2-C(E)$.

We assign each deficient edge to another edge as follows. By Lemma~\ref{lem:deficient}, a deficient edge has a unique triangular common link and contains exactly one low spoke at its apex. There are either one or two low spokes at that apex: three would give the deficient edge two pairs of codegree two. If there are two, assign all deficient wings at this apex to the wing containing both low spokes, called a \emph{type A receiver}. If there is only one, assign them to the wing containing both high spokes, called a \emph{type B receiver}. Each deficient edge is assigned exactly once, to an edge different from itself. For each receiving edge $F$, let $\mathcal R(F)$ be the set of deficient edges assigned to $F$.

\begin{lemma}\label{lem:receiver-payment}
Every receiving edge $F$ satisfies
\begin{equation}\label{eq:receiver-payment}
 C(F)-2\ge\sum_{E\in\mathcal R(F)}\delta(E).
\end{equation}
In particular, every receiving edge is nondeficient.
\end{lemma}

\begin{proof}
Suppose first that $F$ is a type A receiver. Relabel the vertices so that
\[
 E(\CN(a,z))=\bigl\{\{b,c\},\{b,x\},\{c,x\}\bigr\},
 \qquad d_{\Hh}(\{a,b\})=d_{\Hh}(\{a,x\})=2.
\]
The receiver is $F=\{a,b,x\}$. Its possible deficient neighbors at apex $a$ are $\{a,b,c\}$ and $\{a,c,x\}$, and at least one is deficient. Since $N(\{a,b\})=\{c,x\}$ and $N(\{a,x\})=\{b,c\}$, two parts of $K_F$ are singletons with the same original label $c$. The third part consists of tagged copies of the vertices of $N(\{b,x\})\setminus\{a\}$ and contains the copy of $z$, which is adjacent to both singleton copies. If the copy of another vertex $r$ were adjacent to both, then $E(\CN(a,r))=\bigl\{\{b,c\},\{b,x\},\{c,x\}\bigr\}$ as well. The deficient neighbor would have two triangular common links, contrary to Lemma~\ref{lem:deficient}. Thus the copy of $z$ is their unique common neighbor. The singleton copies are nonadjacent, since their labels agree.

Consequently, $K_F$ is obtained from $K_{1,2}$, with one vertex in each part, by adding leaves or isolated vertices to the part containing the copy of $z$. Lemma~\ref{lem:isolates-leaves} gives
\[
 C(F)=\Phi(K_F)\ge\Phi(K_{1,2})=\frac{13}{6}.
\]
Moreover, the copy of $z$ is the only switch, so $F$ has exactly one triangular common link. Hence every request received by $F$ comes from this root. A wing of this root belongs to only one of its two apices. At that apex, the conditions of one and two low spokes are mutually exclusive, so $F$ cannot also be a type B receiver. It follows that $\mathcal R(F)$ consists of at most the two deficient neighbors already identified. Each has deficit smaller than $1/12$ by Lemma~\ref{lem:deficient}, and the surplus of $F$ is at least $1/6$. This proves~\eqref{eq:receiver-payment} for every type A receiver.

Now let $G=\{a,b,c\}$ be a type B receiver. Relabel the vertices so that one of its requests comes from
\begin{equation}\label{eq:type-B-root}
\begin{gathered}
 E(\CN(c,r))=\bigl\{\{a,b\},\{a,u\},\{b,u\}\bigr\},\\
 N(\{c,u\})=\{a,b\},\qquad
 d_{\Hh}(\{c,a\}),d_{\Hh}(\{c,b\})\ge3.
\end{gathered}
\end{equation}
The possible deficient neighbors are $\{a,c,u\}$ and $\{b,c,u\}$. In the parts $X,Y,Z$ of $K_G$, the switch $r_X$ has targets $u_Y,u_Z$. Both targets are ordinary and have neighbors only in $X$. Indeed, an edge from $u_Y$ to some $t_Z$ would require $\{c,u,t\}\in E(\Hh)$, where $t\notin\{a,b\}$, contrary to~\eqref{eq:type-B-root}. The argument for $u_Z$ is the same.

Suppose that $E=\{a,c,u\}$ is deficient. Its singleton part is the tagged copy of $N(\{c,u\})\setminus\{a\}=\{b\}$, its switch has label $r$ and lies in the tagged copy of $N(\{a,u\})\setminus\{c\}$, and the other target of this switch is the copy of $b$ in the part corresponding to $N(\{a,c\})\setminus\{u\}$. A label $s$ is adjacent to this latter target precisely when
\[
 \{a,u,s\},\{a,b,s\}\in E(\Hh),\qquad
 s\notin\{a,b,c,u\}.
\]
These are exactly the conditions for $s_X$ to be adjacent to $u_Y$ in $K_G$. The two target degrees are therefore equal, and Lemma~\ref{lem:deficient} yields
\begin{equation}\label{eq:target-deficit}
 \delta(\{a,c,u\})<\frac1{12}
 -\frac1{2(\deg_{K_G}(u_Y)+1)}.
\end{equation}
For $\{b,c,u\}$, the corresponding conditions are $\{b,u,s\},\{a,b,s\}\in E(\Hh)$, with the same exclusions, and the analogous bound uses $u_Z$. In particular, a target carrying a request has degree at least six.

Call a switch of $K_G$ a \emph{source} if it corresponds to at least one assigned deficit, and call a target \emph{charged} if it carries the corresponding request. For a fixed source, its two targets correspond to its two possible deficient neighbors, so each target carries at most one request. Each source has at least one charged target, and both its targets are ordinary.

We show that a part containing a source contains no other switch. Suppose that $X$ contains a source $r_X$ with target label $u$ and another switch $s_X$ with target label $v$. By Lemma~\ref{lem:labels}, the sets of original labels satisfy
\[
 \{u,v\}\subseteq\{t:t_Y\in Y,\ t_Z\in Z\}
 \subseteq\{r,u\}\cap\{s,v\}.
\]
If $u=v$, then $\CN(c,r)$ and $\CN(c,s)$ are both the triangle on $\{a,b,u\}$. Either incoming deficient edge would have two triangular common links, a contradiction. If $u\ne v$, the inclusions force $r=v$ and $s=u$. The nine triples
\[
\begin{gathered}
 \{a,b,c\},\ \{a,b,r\},\ \{a,b,u\},\ \{a,c,r\},\ \{a,c,u\},\\
 \{b,c,r\},\ \{b,c,u\},\ \{a,r,u\},\ \{b,r,u\}
\end{gathered}
\]
then belong to $E(\Hh)$. Hence $E(\CN(a,b))=\bigl\{\{c,r\},\{c,u\},\{r,u\}\bigr\}$, and both possible incoming deficient edges are wings of this root as well as of $\{c,r\}$, again a contradiction.

Thus sources lie in distinct parts, and there are at most three. Their target pairs are disjoint as sets of tagged vertices: a target shared by sources in different parts would have neighbors in both parts and hence would be a switch, contrary to its ordinariness. No source is a target, since sources are switches and targets are ordinary. Consequently, each member of $\mathcal R(G)$ corresponds to exactly one charged target of exactly one source.

Let $\mathcal S$ be the set of sources and delete them all from $K_G$. This induced deletion is realized by deleting their completion triples from $\Hh$: these triples contain two vertices of $G$, whereas each auxiliary adjacency is witnessed by a triple containing exactly one. The resulting graph $K'$ is therefore an exact tagged graph in a $\mathcal C_4^3$-free subhypergraph retaining $G$.

Every part containing a source has original size at least six, since it contains all neighbors of a charged target. For any one request, the other two parts correspond to its two high spokes and have size at least two. Thus every original part has size at least two. A part containing a source loses one vertex and retains at least five; every other part is unchanged. In particular, all parts of $K'$ have size at least two.

Each target loses exactly one neighbor, and the other remaining degrees are unchanged. Let $J$ be the number of targets of original degree one. These targets are uncharged and distinct, and they become isolated in $K'$. By Lemma~\ref{lem:isolate-surplus},
\begin{equation}\label{eq:deleted-source-charge}
 \Phi(K')\ge2+\frac J6.
\end{equation}
Other isolated vertices can only strengthen this bound.

For each $s\in\mathcal S$, let $h_s$ be the original size of its part and let $d_{s,1},d_{s,2}$ be the original degrees of its targets. Restoring $s$ changes the charge by
\begin{equation}\label{eq:source-restoration}
 \Delta_s=\frac16-\frac1{h_s(h_s+1)}
 -\sum_{j=1}^2\frac1{2d_{s,j}(d_{s,j}+1)}.
\end{equation}
The terms are its switch contribution, the loss in its part contribution, and the losses at its targets. These changes add without interaction: the sources occupy distinct parts, their target sets are disjoint, and no source is a target. Hence the values in~\eqref{eq:source-restoration} are the same regardless of the order of restoration, and
\[
 \Phi(K_G)=\Phi(K')+\sum_{s\in\mathcal S}\Delta_s.
\]

For $j\in\{1,2\}$, let $\delta_{s,j}$ be the deficit assigned to target $j$ of $s$, or zero if that target is uncharged. If it is charged, then~\eqref{eq:target-deficit} and $6\le d_{s,j}\le h_s$ give
\[
 \frac1{2d_{s,j}(d_{s,j}+1)}+\delta_{s,j}
 <\frac1{12}-\frac{d_{s,j}-1}{2d_{s,j}(d_{s,j}+1)}
 \le\frac1{12}-\frac1{h_s(h_s+1)}.
\]
The last step follows from
\[
 \frac{d_{s,j}-1}{2d_{s,j}(d_{s,j}+1)}
 \ge\frac1{d_{s,j}(d_{s,j}+1)}
 \ge\frac1{h_s(h_s+1)}.
\]
For an uncharged target, its loss is $1/4=1/12+1/6$ when $d_{s,j}=1$, and at most $1/12$ when $d_{s,j}\ge2$. Let $j_s$ be the number of degree-one targets of $s$. Since at least one target is charged, summing the two bounds yields
\[
 \sum_{j=1}^2\left(
 \frac1{2d_{s,j}(d_{s,j}+1)}+\delta_{s,j}\right)
 \le\frac16+\frac{j_s}{6}-\frac1{h_s(h_s+1)}.
\]
Together with~\eqref{eq:source-restoration}, this gives
$\Delta_s\ge\delta_{s,1}+\delta_{s,2}-j_s/6$.
The disjointness of the target sets and the correspondence between requests and charged targets give
\[
 \sum_{s\in\mathcal S}j_s=J,
 \qquad
 \sum_{s\in\mathcal S}\sum_{j=1}^2\delta_{s,j}
 =\sum_{E\in\mathcal R(G)}\delta(E).
\]
Therefore~\eqref{eq:deleted-source-charge} gives
\[
\begin{aligned}
 C(G)-2
 &=\Phi(K')-2+\sum_{s\in\mathcal S}\Delta_s\\
 &\ge\frac J6+\sum_{s\in\mathcal S}
   \left(\delta_{s,1}+\delta_{s,2}-\frac{j_s}{6}\right)\\
 &=\sum_{E\in\mathcal R(G)}\delta(E).
\end{aligned}
\]
This proves~\eqref{eq:receiver-payment} for every type B receiver and completes the proof.
\end{proof}

\begin{proof}[Proof of Theorem~\ref{thm:main}]
Each deficient edge is assigned to exactly one receiving edge. By Lemma~\ref{lem:receiver-payment}, every receiver is nondeficient and its surplus covers the sum of all deficits assigned to it. Thus no receiving edge sends a deficit, and no deficit or surplus is counted twice. Every other nondeficient edge has nonnegative surplus. Summing gives
\[
 \sum_{E\in E(\Hh)}(C(E)-2)\ge0.
\]
The charge identity~\eqref{eq:charge-sum} now yields
$2|E(\Hh)|\le|\partial\Hh|+|D(\Hh)|$, as required.
\end{proof}


\begin{thebibliography}{99}
		
		\bibitem{brown}
		W.~G. Brown,
		On graphs that do not contain a Thomsen graph,
		\emph{Canad. Math. Bull.} \textbf{9} (1966), 281--285.
		\href{https://doi.org/10.4153/CMB-1966-036-2}{doi:10.4153/CMB-1966-036-2}.
		
		\bibitem{Erdos3}
		P. Erd\H{o}s,
		On sequences of integers no one of which divides the product of two others and on some related problems,
		\emph{Mitt. Forsch.-Inst. Math. Mech. Univ. Tomsk} \textbf{2} (1938), 74--82.
		\url{https://www.renyi.hu/~p_erdos/1938-07.pdf}.
		
		\bibitem{Erdos}
		P. Erd\H{o}s,
		Problems and results in combinatorial analysis,
		\emph{Congr. Numer.} \textbf{19} (1977), 3--12.
		\url{https://www.renyi.hu/~p_erdos/1977-24.pdf}.
		
		\bibitem{Erdos2}
		P. Erd\H{o}s,
		Some recent problems and results in graph theory,
		\emph{Discrete Math.} \textbf{164} (1997), 81--85.
		\href{https://doi.org/10.1016/S0012-365X(96)00044-1}{doi:10.1016/S0012-365X(96)00044-1}.
		
		\bibitem{ERS}
		P. Erd\H{o}s, A. R\'enyi and V.~T. S\'os,
		On a problem of graph theory,
		\emph{Studia Sci. Math. Hungar.} \textbf{1} (1966), 215--235.
		\url{https://real.mtak.hu/110437/1/1966-06.pdf}.
		
		\bibitem{FurediC4}
		Z. F\"uredi,
		Graphs without quadrilaterals,
		\emph{J. Combin. Theory Ser. B} \textbf{34} (1983), 187--190.
		\href{https://doi.org/10.1016/0095-8956(83)90018-7}{doi:10.1016/0095-8956(83)90018-7}.
		
		\bibitem{Furedi}
		Z. F\"uredi,
		Hypergraphs in which all disjoint pairs have distinct unions,
		\emph{Combinatorica} \textbf{4} (1984), 161--168.
		\href{https://doi.org/10.1007/BF02579216}{doi:10.1007/BF02579216}.
		
		\bibitem{Hanani}
		H. Hanani,
		On balanced incomplete block designs with blocks having five elements,
		\emph{J. Combin. Theory Ser. A} \textbf{12} (1972), 184--201.
		\href{https://doi.org/10.1016/0097-3165(72)90035-0}{doi:10.1016/0097-3165(72)90035-0}.
		
		\bibitem{HMY}
		H. Huang, J. Ma and T. Yang,
		Extremal hypergraphs without generalized $4$-cycles,
		arXiv:2609.37744 (2026).
		\url{https://arxiv.org/abs/2609.37744}.
		
		\bibitem{MY}
		J. Ma and T. Yang,
		Upper bounds on the extremal number of the $4$-cycle,
		\emph{Bull. Lond. Math. Soc.} \textbf{55} (2023), 1655--1667.
		\href{https://doi.org/10.1112/blms.12810}{doi:10.1112/blms.12810}.
		
		\bibitem{Marchant}
		E. Marchant,
		The largest minimum codegree of a $3$-graph without a generalized $4$-cycle,
		\emph{Combin. Probab. Comput.} \textbf{22} (2013), 112--117.
		\href{https://doi.org/10.1017/S0963548312000430}{doi:10.1017/S0963548312000430}.
		
		\bibitem{MubayiExact}
		D. Mubayi,
		Some exact results and new asymptotics for hypergraph Tur\'an numbers,
		\emph{Combin. Probab. Comput.} \textbf{11} (2002), 299--309.
		\href{https://doi.org/10.1017/S0963548301005028}{doi:10.1017/S0963548301005028}.

			\bibitem{MVbipartite}
		D. Mubayi and J. Verstra\"ete,
		A hypergraph extension of the bipartite Tur\'an problem,
		\emph{J. Combin. Theory Ser. A} \textbf{106} (2004), 237--253.
		\url{https://www.math.cmu.edu/~mubayi/papers/hyperbtpf.pdf}.
		
		\bibitem{MVminimal}
		D. Mubayi and J. Verstra\"ete,
		Minimal paths and cycles in set systems,
		\emph{European J. Combin.} \textbf{28} (2007), 1681--1693.
		\href{https://doi.org/10.1016/j.ejc.2006.07.001}{doi:10.1016/j.ejc.2006.07.001}.
		
		\bibitem{PV}
		O. Pikhurko and J. Verstra\"ete,
		The maximum size of hypergraphs without generalized $4$-cycles,
		\emph{J. Combin. Theory Ser. A} \textbf{116} (2009), 637--649.
		\href{https://doi.org/10.1016/j.jcta.2008.09.002}{doi:10.1016/j.jcta.2008.09.002}.
		
	\end{thebibliography}
\end{document}